\documentclass[12pt]{article}
\usepackage{latexsym,amsmath,amsopn,amssymb,amsthm,amsfonts}

\begin{document}

\title {Construction of some types wavelets with coefficient of
scaling $N$}
\author{N.~K.~Smolentsev, P.~N.~Podkur
\thanks{Kemerovo State University, Kemerovo, 650043, RUSSIA.
E-mail smolen@kuzbass.net} } \maketitle

\begin{abstract}
In this paper it is shown, that $B$-splines are scaling functions
for any natural $N$. Construction of Haar wavelets,
Kotelnikov-Shannon wavelets, nonorthogonal wavelets based on
$B$-splines is given.
\end{abstract}

{\bf 1. Scaling functions.} Let $N>1$ -- an integer, $\mathbb{Z}$
-- the set of all integers and $L^2(\mathbb{R})$ -- Hilbert space
of square integrable functions.

{\bf Definition 1.} {\it Function  $\varphi(x) \in
L^2(\mathbb{R})$ is called $N$-scaling, if it can be represented
as
$$
\varphi(x)=\sqrt{N}\sum_{n\in \mathbb{Z}} h_n \varphi(Nx-n),
\eqno(1)
$$
where coefficients $h_n$,\ $n\in \mathbb{Z}$ satisfy to a
condition $\sum_{n}|h_n|^2<\infty$. The relationship (1) is called
the $N$-scale equation (refinement equation). The set $\{h_n\}$ of
coefficients of expansion in the equation (1) is called the
scaling filter.}

{\bf Note 1.} If $N$-scaling function $\varphi(x)$ has the compact
support of length $L$, then the sum in equation (1) is finite and
consist at most $L(N-1)+1$ components.

The Fourier transform of $N$-scale equation is
$$
\widehat{\varphi}(\omega)=H_0\left( \frac{\omega}{N}\right)
\widehat{\varphi}\left(\frac{\omega}{N} \right) , \eqno(2)
$$
where
$$
H_0(\omega)=\frac {1}{\sqrt{N}} \sum_{n\in \mathbb{Z}} h_n
e^{-in\omega}. \eqno(3)
$$
Function $H_0(\omega)$ we shall name {\it frequency function} of
scaling function $\varphi(x)$.

Let's give examples of $N$-scaling functions.

{\bf 1.1.  Haar function.} It is defined as characteristic
function of an interval $[0,1)$,
$$
\varphi_0(x)=
\left\{%
\begin{array}{l}
  1, \qquad x\in [0,1), \\
  0, \qquad x\notin [0,1) \\
\end{array}%
\right. .
$$
Fourier transform of the function is
$$
\widehat{\varphi}_0(\omega)=e^{-i\omega/2} \frac {\sin
\omega/2}{\omega/2}.  \eqno(4)
$$
It is easy to see, that for any natural $N$
$$
\varphi_0(x)=\varphi_0(Nx)+\varphi_0(Nx-1)+\dots
\varphi_0(Nx-N+1)= \sqrt{N}\sum_{n=0}^{N-1}h_n \varphi_0(Nx-n),
\quad h_n=\frac {1}{\sqrt{N}}.
$$
Frequency function is
$$
H_0(\omega)= \frac {1}{N}\left(1+ e^{-i\omega}+
e^{-i2\omega}+\dots  e^{-i(N-1)\omega} \right) = \frac
{1}{N}e^{-i(N-1)\omega/2} \frac {\sin N\omega/2}{\sin \omega/2} .
$$

{\bf 1.2. Kotelnikov-Shannon scaling function.} Let's show, that
function $\varphi(x)= \frac {\sin \pi x}{\pi x}$ is $N$-scaling
for any natural $N>1$.

{\bf Theorem 1}.{ \it The function  $\varphi(x)= \frac {\sin \pi
x}{\pi x}$ is $N$-scale function for any natural $N>1$. The
refinement equation for the function has form:
$$
\varphi(x)=\sum_{n\in \bf{Z}}\varphi(n/N)\frac {\sin\pi (Nx-n)}{\pi (Nx-n)}=
\sqrt{N}\sum_{n\in \bf{Z}}h_n \frac {\sin\pi (Nx-n)}{\pi (Nx-n)},
\eqno(5)
$$
where
$$
h_n=\frac{\sqrt{N}}{\pi n}\sin\left(\frac {\pi n}{N}\right).
$$}

{\it Proof}. We shall remind, that function  $f(x) \in
L^2(\mathbb{R})$ is called {\it band limited} function \cite{Db},
if its Fourier transform $\widehat{f}(\omega)$ is equal to zero
outside of some frequency band $[-\Omega, \Omega]$. The Fourier
transform $\widehat{\varphi}(\omega)$ of function $\varphi(x)=
\frac {\sin \pi x}{\pi x}$ has the support on interval
$[-\pi,\pi]$ and $\widehat{\varphi}(\omega)\equiv 1$ on this
interval $[-\pi,\pi]$. Then the support belongs to the greater
interval $[-\pi N,\pi N]$. From Kotelnikov theorem \cite{Sm}, band
limited function $f(x)$ of a band $[-\Omega,\Omega]$ can be
represented as
$$
f(x)=\sum_{n\in \bf{Z}} f \left(\frac {\pi}{\Omega}n \right) \frac
{\sin (\Omega x-\pi n)}{\Omega x-\pi n}.
$$
In particular, for our function  $\varphi (x)= \frac {\sin \pi
x}{\pi x}$ if $\Omega=\pi N$, then
$$
\varphi(x)=\sum_{n\in \bf{Z}}\varphi \left( \frac {\pi}{\pi N}n
\right) \frac {\sin (\pi Nx-\pi n)}{ (\pi Nx-\pi n)}= \sum_{n\in
\bf{Z}}\varphi \left( \frac {n}{N}\right) \frac {\sin \pi( Nx-n)}{
\pi( Nx-n)}= \sqrt{N}\sum_{n\in \bf{Z}}h_n \varphi(Nx-n).
$$

{\bf 1.3. $B$-spline scaling functions.} Well-known \cite{Db},
\cite{Sm}, that $B$-splines are 2-scaling functions. Let's prove,
that $B$-splines are $N$-scaling for any natural $N>1$. The
$B$-spline of order zero is Haar function $\varphi_0(x)$
(characteristic function of an interval $[0,1)$). $B$-spline
$\varphi_1(x)$ of order one with support on interval $[0,2]$ can
be presented \cite{Sm} as convolution $\varphi_0(x)$ with itself,\
$\varphi_1(x)=\varphi_0(x)\ast\varphi_0(x)$. The higher-order
$B$-splines are defined inductively,\
$\varphi_n(x)=\varphi_{n-1}(x)\ast\varphi_0(x)$. Further shall
consider, that support of $B$-spline $\varphi_n(x)$ is on the
interval $[-(n +1)/2, (n +1)/2]$ -- in case of odd $n$, and in
case of even $n$ -- on the interval $[-n/2, n/2+1]$. Then the
Fourier transform of $B$-spline $\varphi_n(x)$ has form \cite{Db},
\cite{Sm}:
$$
\widehat{\varphi}_n(\omega)=e^{-iK\omega/2}\left(\frac {\sin
\omega/2}{\omega/2}\right)^{n+1}. \eqno(6)
$$
where $K=0$ -- in case of odd $n$ and $K=1$ -- in case of even $n$.

{ \bf Theorem 2}. {\it $B$-spline $\varphi_n(x)$ of order $n$ is
$N$-scaling function for any natural $N>1$ and any natural $n$.
The right part of scaling equation $\varphi(x)=\sqrt{N}\sum_{k}
h_n \varphi(Nx-k)$ is the finite sum. Filter coefficients $h_k$
are find from the formula
$$
h_k=\frac{1}{\sqrt{N}
N^n}\sum_{\alpha}\frac{(n+1)!}{\alpha_0!\alpha_1!\dots
\alpha_{N-1}!},\eqno (7)
$$
where in case of odd $n$ summation is made
on all multi-indexes $\alpha=(\alpha_0,\alpha_1,\dots \alpha_{N-1})$,
satisfying two conditions:
$$|\alpha|=\alpha_0 +\alpha_1 +\dots + \alpha_{N-1}= n+1
$$
$$\alpha_1 +2\alpha_2+\dots +
(N-1)\alpha_{N-1}-(N-1)(n+1)/2 = k,
$$
the index $k$ varies from $k =
-(N-1)(n+1)/2$ up to $k = (N-1)(n+1)/2$.

In case of even $n$, filter coefficients $h_k$ are find from the
same formula (7), where summation is made on all multi-indexes
$\alpha=(\alpha_0,\alpha_1,\dots \alpha_{N-1})$, satisfying two
conditions:
$$|\alpha|=\alpha_0 +\alpha_1 +\dots + \alpha_{N-1}= n+1
$$
$$\alpha_1 +2\alpha_2+\dots + (N-1)\alpha_{N-1}-(N-1)n/2 =k,
$$
the index $k$ varies from $k = -(N-1)n/2$ up to $k = (N-1)n/2+1$. }

{\it Proof}. First shall consider a case odd $n$. Then the support
of $\varphi_n(x)$ is on the interval $[-(n+1)/2,  (n  +1)/2]$ and
the Fourier transform looks like (6), where $K=0$. Let's
substitute (6) in scaling equation (2) and shall find frequency
function $H_0(\omega)$,
$$
H_0(\omega)=\frac{\widehat{\varphi}_n(N\omega)}{\widehat{\varphi}_n(\omega)}=
\left(\frac {\sin(N\omega/2)}{N\omega/2}\right)^{n+1}
\left(\frac{\omega/2}{\sin(\omega/2)}\right)^{n+1}= \left(\frac 1
N \frac{\sin(N\omega/2)}{\sin(\omega/2)}\right)^{n+1}=
$$
$$
=\left(\frac 1 N \frac {e^{iN\omega/2}-e^{-iN\omega/2}}
{e^{i\omega/2}-e^{-i\omega/2}}\right)^{n+1}= \left(\frac 1 N
e^{i(N-1)\omega/2}\right)^{n+1}\left(\frac
{1-e^{-iN\omega}}{1-e^{-i\omega}}\right)^{n+1}=
$$
$$
=\frac 1 {N^{n+1}} e^{i(N-1)(n+1)\omega/2} \left(1+e^{-i\omega}+
\dots +e^{-i(N-1)\omega}\right)^{n+1}.
$$
Since  $n$ is odd, the multiplier $e^{i\omega(N-1)(n+1)/2}$ is an
integer degree of $e^{i\omega}$. Expanding a degree
$\left(1+e^{-i\omega}+ \dots +e^{-i(N-1)\omega}\right)^{n+1}$, we
obtain frequency function $H_0(\omega)$ as the finite sum on
degrees of a variable $z=e^{-i\omega}$, i.e.
$H_0(\omega)=\sqrt{N^{-1}} \sum_k h_k e^{-ik\omega}$. Factors of
expansion of function $H_0(\omega)$ give us required filter
coefficients $h_k$ of a scaling equation. Thus summation is made
from $k =  -(N-1)(n+1)/2$ up to $k =(N-1)(n+1)/2$. We shall
calculate these coefficients. For this purpose we use known analog
of a Newton binomial formula
$$
(z_0+z_1+\dots+z_p)^m=\sum_{|\alpha|=m}\frac{m!}
{\alpha_0!\dots\alpha_p!}z_0^{\alpha_0}\dots z_p^{\alpha_p}.
$$
Then
$$
\left(1+e^{-i\omega}+\dots+e^{-i(N -1)\omega}\right)^{n+1}=
\sum_{k=0}^{(N-1)(n+1)}\sum_{\scriptsize \begin{array}{c}
  |\alpha|= n+1, \\
  \alpha_1{+}2\alpha_2{+}\dots{+}(N{-}1)\alpha_{N-1}=k
\end{array}}\! \frac{(n+1)!}{\alpha_0!\alpha_1!\dots\alpha_{N-1}!}e^{-ik\omega},
$$
$$
H_0(\omega)=\frac 1 {N^{n+1}} e^{i(N-1)(n+1)\omega/2}
\left(1+e^{-i\omega}+ \dots +e^{-i(N-1)\omega}\right)^{n+1}=
$$
$$
=\frac 1{N^{n+1}}\sum_{k=-(N-1)(n+1)/2}^{(N-1)(n+1)/2}
\sum_{\scriptsize \begin{array}{c}
|\alpha|=n+1, \\
\alpha_1+2\alpha_2+\dots+(N-1)\alpha_{N-1}=
k+(N-1)(n+1)/2\end{array}} \frac{(n+1)!}
{\alpha_0!\alpha_1!\dots\alpha_{N-1}!}e^{-ik\omega}.
$$
From last expression the statement of the theorem follows. In the
second case, when $n$ is even, the proof is made similarly. The
only difference that support $\varphi_n(x)$ is on interval $[-n/2,
n/2 +1]$ and Fourier transform of function $\varphi_n(x)$ has a
phase factor $e^{-i\omega/2}$ according to the formula (6).

{\bf Example 1.3.1.} $B$-spline $\varphi_1(x)$ of order 1,
$$
\varphi_1(x)=
\left\{%
\begin{array}{l}
  1-|x|, \qquad x \in [-1,1], \\
  0, \qquad \qquad x \notin [-1,1]\\
\end{array}%
\right. .
$$
Fourier transform of the function is
$$
\widehat{\varphi}_1(\omega)=\left(\frac{\sin
\omega/2}{\omega/2}\right)^2=2 \frac{1-\cos \omega}{\omega^2}.
$$
It is easy to see, that at $N=3$,
$$
\varphi_1(x)=\frac 1 3 \varphi(3x+2)+\frac 2 3 \varphi(3x+1)+
\varphi(3x)+\frac 2 3 \varphi(3x-1)+\frac 1 3 \varphi(3x-2),
$$
$$
H_0(\omega)=\frac 1 {\sqrt{3}}\frac 1 {3\sqrt{3}}\left(
e^{i2\omega}+2e^{i\omega}+3+2e^{-i\omega}+e^{-i2\omega}\right)=
\frac 1 9 \left(3+4\cos\omega+2\cos2\omega\right).
$$
At any natural $N>1$ the following $N$-scaling equation has form
$$
\varphi_1(x)=\sum_{n=-N+1}^{N-1}\frac{N-|n|}{N}\varphi(Nx-n).
$$
Notes, that the family $\{\varphi(x-n),\  n\in \mathbb{Z}\}$ form
basis of space $V_0\subset L^2(\mathbb{R})$, but it is not
orthogonal.

{\bf Example 1.3.2.} Let's consider $B$-spline $\varphi_2(x)$ of order 2,
$$
\varphi_2(x)=
\left\{%
\begin{array}{ll}
  (x+1)^2/2, & x \in [-1,0], \\
  3/4-(x-1/2)^2, & x \in [0,1],\\
  (x-2)^2/2, & x \in [1,2],\\
  0, & x  \notin [-1,2]\ .
\end{array}%
\right.
$$
Fourier transform of the function is
$$
\widehat{\varphi}_2(\omega)=e^{-i\omega/2}\left(\frac {\sin
\omega/2}{\omega/2}\right)^3.
$$
For $N =3$ filter coefficients $h_k$ of the scaling filter is find
from the formula (7). Then,
$$
h_{-2}=\frac{1}{9\sqrt{3}},\ h_{-1}=\frac{1}{3\sqrt{3}}, \
h_{0}=\frac{2}{3\sqrt{3}}, \ h_{1}=\frac{7}{9\sqrt{3}}, \
h_{2}=\frac{2}{3\sqrt{3}}, \ h_{3}=\frac{1}{3\sqrt{3}}, \
h_{4}=\frac{1}{9\sqrt{3}}.
$$
The scaling relation is
$$
\begin{array}{l}
\varphi_2(x)=\frac 19 \varphi_2(3x+2)+ \frac 13 \varphi_2(3x+1)+
\frac 23 \varphi_2(3x)+ \\
\qquad \qquad\qquad +\frac 79 \varphi_2(3x-1)+ \frac 23
\varphi_2(3x-2)+ \frac 13 \varphi_2(3x-3)+ \frac 19
\varphi_2(3x-4).
\end{array}
$$

{\bf 2. $N$-multiresolution analysis}. Here we will define a basic constructions of
multiresolution analysis for a case of arbitrary natural coefficient
of scaling $N> 1$, for more details about it see  \cite{Db}, \cite{BJ}.

Orthogonal $N$-multiresolution analysis of $L^2(\mathbb{R})$ is a
sequence of closed nested (telescoping) subspaces $\dots \subset
V_{-1}\subset V_{0}\subset V_{1}\subset V_{2} \subset\dots$ which
has properties:

1) $\cup_{j \in \mathbb{Z}}V_j$ is dense $L^2(\mathbb{R})$.

2) $\cap_{j \in
\mathbb{Z}}V_j=\{0\}$.

3) $f(x) \in V_j \quad \Leftrightarrow \quad f(Nx) \in V_{j+1}$.

4) There is a function $\varphi(x) \in V_0$,  called scaling function,
that functions $\varphi_{0,n}(x)=\varphi(x-n), \ n
\in \mathbb{Z}$ form an orthonormal basis of space
$V_0$.

From properties 3 and 4 implies, that functions
$$
\varphi_{j,n}=\sqrt{N^j}\varphi(N^jx-n), \quad n \in \mathbb{Z}
\eqno(8)
$$
form an orthonormal basis of space $V_j$ for any $j$.

From $\varphi(x) \in V_0 \subset V_1$ follows that function $\varphi(x)$
is decomposed on basis of space $V_1$, i.e. it is $N$-scaling function,
$$
\varphi(x)=\sqrt{N}\sum_n h_n \varphi(Nx-n).
$$
From the refinement equation follows expansion of base function (8)
for any $j,k \in \mathbb{Z}$:
$$
\varphi_{j-1,k}(x)=\sum_n h_{n-Nk}\varphi_{j,n}(x)= \sum_n
h_n\varphi_{j,n+Nk}(x). \eqno(9)
$$
Let's remind \cite{Db}, that integer translates $\varphi_n(x)=\varphi(x-n)$
form orthonormal basis of subspace $V_0 \subset L^2(\mathbb{R})$,
if and only if
$$
\sum_{n \in {\bf Z}}|\widehat{\varphi}(\omega+2\pi n)|^2 = 1 \quad
\text{\it a.e }.
$$
Hence \cite{Db}, \cite{BJ}, if translates $\varphi(x-n)$ form
orthonormal basis of subspace $V_0$, then frequency function
$H_0(\omega)$ satisfy the following property
$$
\sum_{m=0}^{N-1}|H_0(\omega+2\pi m/N)|^2 = 1 \quad \text{\it constancy a. e.}
$$

{\bf 2.1.  Wavelets.} Well-known \cite{Db}, \cite{BJ}, that in
case $N=2$ in orthogonal multiresolution analysis
$L^2(\mathbb{R})$ each subspace $V_j$ of $V_{j+1}$ has orthogonal
complement $W_j$ and $V_{j+1}=V_j\oplus W_j$. There is a function
$\psi(x)$, called mother wavelet, that the set of its translates
$\psi(x-n)$ forms an orthonormal basis of space $W_0$ and
functions $\psi_{j,n}(x)= \sqrt{2^j}\psi(2^jx-n), \ n \in
\mathbb{Z}$ forms an orthonormal basis of space $W_j$ for any
$j\in \mathbb{Z}$.

In case of arbitrary integer coefficient of scaling $N>1$ situation some differs.
To scaling function $\varphi(x)$ corresponds $N-1$ wavelets
$\psi^1(x),\dots ,
\psi^{N-1}(x)$\ \cite{Db}, \cite{BJ},
$$
\psi^k(x)=\sqrt{N}\sum_n g_n^k\varphi(Nx-n), \quad k=1,2,\dots ,
N-1.
$$
In frequency space we have:
$$
\widehat{\psi}^k(\omega)=G_k\left(\frac {\omega}{N}\right)
\widehat{\varphi}\left(\frac {\omega}{N}\right), \quad k=1,2,\dots
, N-1,
$$
Where $G_k(\omega)=\sqrt{N^{-1}}\sum_n g_n^ke^{-in\omega}$ there
are frequency functions corresponding to wavelets $\psi^k(x)$.
Then the following matrix is unitary \cite{Db}, \cite{BJ}:
$$
\left( \begin{array}{cccc}
  H_0(\omega) & H_0(\omega+2\pi/N) & \ldots & H_0(\omega+2\pi(N-1)/N) \\
  G_1(\omega) & G_1(\omega+2\pi/N) & \ldots & G_1(\omega+2\pi(N-1)/N) \\
  \hdotsfor[2.5]{4} \\
  G_{N-1}(\omega) & G_{N-1}(\omega+2\pi/N) & \ldots & G_{N-1}(\omega+2\pi(N-1)/N) \\
  \end{array} \right) \eqno(10)
$$
Then we get $N-1$ spaces of wavelets: $W_j^k$, $k=1,2,\dots , N-1$
and the orthogonal expansion for any $j$
$$
V_{j-1}=V_j\oplus W_j^1\oplus \dots W_j^{N-1}.
$$
Functions
$$
\psi_{j,n}^k(x)=\sqrt{N^j}\psi^k(x)(N^jx-n),\quad n\in \mathbb{Z}
.
$$
form an orthonormal basis of spaces of wavelets $W_j^k$,
$i=1,2,\dots , N-1$.

{\bf 2.2. Examples of wavelets.} We shall present two examples of
$N$-scaling functions and corresponding wavelets.

{\bf 2.2.1. Haar wavelets.} For Haar function $\varphi(x)$\
(characteristic function of an interval $[0,1)$), $N$-scaling
equation is indicated in 1.1. Nonzero coefficients of the filter
$h_n$ there are: $h_0=h_1 =\dots = h_{N-1} =1/\sqrt{N}$.

We shall find wavelet-functions $\psi^1(x),\dots, \psi^{N-1}(x)$.
Let  $\psi^i(x)\in W_0^i$. As $W_0^i\subset V_1$, then $\psi^i(x)$ it is decomposed on basis
$\varphi_{1,n}(x)$ of $V_1$: $\psi^i=\sum_ k c_n^i\varphi_{1,k}$. As $\psi^i(x)\perp V_0$,
then for any $n$ we get: $\left(\psi^i, \varphi_{0,n}\right)= 0$.
Space $V_0$ is nested to $V_1$, hence functions  $\varphi_{0,n}(x)$ also are decomposed
on basis $\varphi_{1,n}(x)$ of $V_1$ according to the formula (9):
$\varphi_{0,n}=\sum_p h_k\varphi_{1,Nn+p}$.
More in detail
$$
\varphi_{0,n} = \frac {1}{\sqrt{N}}\varphi_{1,Nn}+ \frac
{1}{\sqrt{N}}\varphi_{1,Nn+1}+ \dots+ \frac
{1}{\sqrt{N}}\varphi_{1,Nn+N-1}.
$$
Then the condition of orthogonality $\psi^i$ to $V_0$ becomes
$$
\left(\psi^i, \varphi_{0,n}\right) = \left(\sum_k
c_k^i\varphi_{1,k}, \frac
{1}{\sqrt{N}}\sum_{p=0}^{N-1}\varphi_{1,Nn+p}\right)=0.
$$
As the family $\{\varphi_{1,n}(x)\}$ form orthonormal basis, then
from last equality we shall get
$$
c_{Nn}^i+ c_{Nn+1}^i+ \dots+ c_{Nn+N-1}^i=0, \quad \forall\ n,
\quad i=1,2,\dots , N-1.
$$
From condition of orthogonality of wavelets, $\left(\psi^i, \psi^j
\right) =\delta_{ij}$  follows equalities $\sum_n c_{n}^i c_{n}^j=\delta_{ij}$.
We obtain the following set of equations:
$$
\left\{\begin{array}{l} c_{Nn}^i+ c_{Nn+1}^i+ \dots+
c_{Nn+N-1}^i=0, \quad \forall\ n, \quad i=1,2,\dots , N-1,\\
\sum_n c_{n}^i c_{n}^j=\delta_{ij}, \quad i,j =1,2,\dots , N-1 .
\end{array} \right. \eqno(11)
$$
The system has a variety of solutions. The elementary solution corresponding
to a minimum set of nonzero coefficients is ambiguously. A solution will be any set of $N-1$
orthogonal vectors in a plane in $\mathbb{R}^N$, specified by the equation $x_0+x_1+\dots +x_{N-1}=0$

To solving of system we shall construct an orthogonal matrix of
dimension $N$, first row of which is a vector $h=(h_0,h_1,\dots ,
h_{N-1})=(1,1,\dots , 1)/\sqrt{N}$. It is possible to take, for
example, the following $N$-dimensional matrix, $N$,
$$
\left(%
\begin{array}{cccccc}
  1 & 1 & 1 & \ldots & 1 & 1 \\
  1 & -1 & 0 & \ldots & 0 & 0 \\
  1 & 1 & -2 & \ldots & 0 & 0 \\
  \hdotsfor[2.5]{6} \\
  1 & 1 & 1 & \ldots & 1 & -N+1 \\
\end{array}%
\right)
$$
and then to make a normalization of these rows.

For example, if $N=3$ we can to select the elementary solution:
$c_0^1=1,\ c_1^1=-1,\ c_2^1=0$ and $c_0^2=1,\ c_1^2=1,\ c_2^2=-2$,
remaining coefficients we shall accept equal to zero. Normalize these vectors:
$$
c_0^1=\frac{1}{\sqrt{2}},\ c_1^1=-\frac{1}{\sqrt{2}},\
c_2^1=0\quad \text{ and } \quad  c_0^2=\frac{1}{\sqrt{6}},\
c_1^2=\frac{1}{\sqrt{6}},\ c_2^2= -\frac{\sqrt{2}}{\sqrt{3}}.
$$
Then
$$
\psi^1(x)=\frac{1}{\sqrt{2}} \varphi_{1,0}(x)-\frac{1}{\sqrt{2}}
\varphi_{1,1}(x) =\frac{\sqrt{3}}{\sqrt{2}}
\varphi(3x)-\frac{\sqrt{3}}{\sqrt{2}} \varphi(3x-1),
$$
and
$$
\psi^2(x)=\frac{1}{\sqrt{2}}\varphi(3x)+ \frac{1}{\sqrt{2}}
\varphi(3x-1) -\sqrt{2}\varphi(3x-2).
$$

{\bf 2.2.2. Kotelnikov-Shannon wavelets.} Let's consider Kotelnikov function
$$
\varphi(x)= \frac {\sin \pi x}{\pi x},\quad
\widehat{\varphi}(\omega)=
\left\{%
\begin{array}{l}
  1, \quad \omega \in [-\pi,\pi], \\
  0, \quad \omega \notin [-\pi,\pi]
\end{array}%
\right.
$$
which, as shown in the theorem 1, is $N$-scaling function. Let's
construct $N$-multiresolution analysis. Space $V_0$, generated by
$\varphi(x-n)$, $n\in \mathbb{Z}$, consist of all square
integrable functions $f(x)$ for which Fourier transform has
support on interval $[-\pi,\pi]$. Space $V_j$ consist of functions
whose Fourier transforms have the support on interval
$[-N^j\pi,N^j\pi]$. Let's find wavelets
$\psi^1(x),\dots,\psi^{N-1}(x)$. First we shall define frequency
function $H_0(\omega)$, and then shall find functions
$G_k(\omega)$ from unitary condition of a matrix (10). Function
$H_0(\omega)$ can be found from the scaling equation,
$\widehat{\varphi}(\omega)=H_0(\omega/N)\widehat{\varphi}(\omega/N)$,
$$
H_0(\omega)=\left\{\begin{array}{l}
  1, \qquad \omega \in [-\pi/N,\pi/N], \\
  0, \qquad \omega \in [-\pi,-\pi/N]\cup [\pi/N,\pi] .
\end{array} \right.
$$
Outside of interval $[-\pi,\pi]$ function $H_0(\omega)$ extends
periodically For $k =1,2,\dots , N-1$ we shall put
$$
G_k(\omega)=\left\{\begin{array}{l}
  1, \quad |\omega| \in [k\pi/N,(k+1)\pi/N], \\
  0, \quad \text{ for other }\ \omega \in [-\pi,\pi].
\end{array} \right.
$$
Wavelet $\psi^k\in W_0^k$ is defined from the formula
$$
\widehat{\psi}^k(\omega)= G_k(\omega/N)\widehat{\varphi}(\omega/N)
=\left\{\begin{array}{l}
  1, \quad |\omega| \in [k\pi,(k+1)\pi], \\
  0, \quad \text{ for other }\ \omega .
\end{array} \right.
$$
Then for $k =1,2,\dots , N-1$ we have,
$$
\psi^k(x)=\frac{1}{2\pi}\int_{-\infty}^{\infty}\widehat{\psi}^k(\omega)
e^{i\omega x}\ d\omega
=\frac{1}{2\pi}\int_{-(k+1)\pi}^{-k\pi}\widehat{\psi}^k(\omega)
e^{i\omega x}\ d\omega
+\frac{1}{2\pi}\int_{k\pi}^{(k+1)\pi}\widehat{\psi}^k(\omega)
e^{i\omega x}\ d\omega =
$$
$$
=\frac{\sin (k+1)\pi x- \sin k\pi x}{\pi x}\ .
$$
This function we shall name {\it Kotelnikov-Shannon wavelets} of
scale $N$.

{\bf 3.  Wavelet-transform.} First we shall assume, that we deal
with the orthonormal wavelets. Then we have scaling function
$\varphi(x)$ and wavelets $\psi^1(x),\dots,\psi^{N-1}(x)$.
Wavelet-expansion of signal $A = \{a_n\}$ is  \cite{Db},
\cite{BJ}:
$$
a_{1,k}=\sum_n \overline{h_n} a_{n+Nk}, \qquad d_{1,k}^i=\sum_n
\overline{g_n^i} a_{n+Nk}, \quad i=1,2,\dots , N-1. \eqno(12)
$$
where $\{h_n,\ n\in \mathbb{Z}\}$ and $\{g_n^i,\ n\in
\mathbb{Z}\}$ there are filters of wavelets $\varphi(x)$ and
$\psi^i(x)$. Let's write last formulas as convolution. For this
purpose we shall enter the following coefficients:
$h_n^*=\overline{h_{-n}}$, $g_{n}^{i*}=\overline{g_{-n}^i}$. Then
$$
a_{1,k}=\sum_n h_n^* a_{Nk-n}, \qquad d_{1,k}^i=\sum_n g_{n}^{i*}
a_{Nk-n}, \quad i=1,2,\dots , N-1. \eqno(13)
$$
Thus, wavelet-expansion is made by the conjugate filters
$\{h_n^*\}$ and $\{g_{n}^{i*}\}$, with next $N$-adic decimation (a
choice only elements with numbers $Nk$). Procedure of expansion
can be repeated, having applied it to set of coefficients
$cA_1=\{a_{1,k}\}$. Reconstruction of array $A = \{a_n\}$ by
wavelet-coefficients is
$$
a_n= \sum_k h_{n-Nk} a_{1,k}+\sum_{i=1}^{N-1}\sum_k g_{n-Nk}^i
d_{1,k}^i . \eqno(14)
$$
Last formula can be written also as convolution, having made
inverse $N$-adic decimation of arrays $cA_1=\{a_{1,k}\}$ and
$cD_1=\{d_{1,k}^1, \dots d_{1,k}^{N-1}\}$,
$$
\widetilde{a}_{1,m}=\left\{\begin{array}{l}
  a_{1,k}, \quad \text{ if } m=Nk, \\
  0, \quad \text{ if } m \text{ not multiply } Nk
\end{array} \right. , \qquad
\widetilde{d}_{1,m}^{\ i} =\left\{\begin{array}{l}
  d_{1,k}^{\ i}, \quad \text{ if } m=Nk, \\
  0, \quad \text{ if } m \text{ not multiply } Nk
\end{array} \right.
$$
Then the formula (14) becomes
$$
a_n= \sum_k h_{n-m} \widetilde{a}_{1,m}+\sum_{i=1}^{N-1}\sum_k
\widetilde{g}_{n-m}^{\ i} d_{1,m}^i . \eqno(15)
$$

{\bf 3.2. Decomposition and reconstruction in not orthogonal
case.} Formulas of wavelet decomposition and reconstruction (12)
-- (15) usable only for an orthogonal case. They are unadaptable,
if expansion $V_j=V_{j-1}\oplus W_{j-1}^1 \oplus \dots \oplus
W_{j-1}^{N-1}$ is not orthogonal and functions
$$
\varphi_{j,n}(x)=\sqrt{N^j}\varphi(Nx-n),\quad
\psi_{j,n}^i(x)=\sqrt{N^j}\psi^i(Nx-n),\quad n\in \mathbb{Z},\
i=1,2,\dots, N-1
$$
do not form the orthonormal systems. Let decomposition of a signal
$\{a_k\}$ is made by some (nonorthogonal) filters $\{h_k\}$ and
$\{g_k^i\}$, $i = 1, 2,\dots , N-1$, by (12),
$$
a_{1,k}=\sum_n h_n a_{Nk-n},\quad d_{1,k}^i=\sum_n g_n^i
a_{Nk-n},\quad i=1,2,\dots, N-1
$$
Let's find other filters $\{\widetilde{h}_k\}$ and
$\{\widetilde{g}_k^i\}$, $i = 1, 2,\dots , N-1$, which lead to
perfect reconstruction of the input signal by (15). It is
convenient to solve the task in terms of formal power series.

Let $X(z)=\sum_n a_n z^n$ is the series corresponding to the
signal $\{a_k\}$. We shall define frequency functions
$H_0(z)=\sum_n h_n z^n$ and $H_i(z)=\sum_n g_n^i z^n$, $i = 1,
2,\dots , N-1$. Then their effect on a signal is defined by
multiplication:
$$
X_0(z)=H_0(z)X(z),\quad X_i(z)=H_i(z)X(z),\quad i=1,2,\dots, N-1 .
$$
At wavelet expansion it is still necessary made sample of elements
with numbers $Nk$. In terms of formal power series this procedure
is reduced to a choice of elements of power series with degrees,
multiple $N$. It can be made as follows. Let $\rho=e^{i2\pi/N}$,
then it is easy to see, that the following sum contains only
degrees, multiple $z^N$,
$$
\frac 1N\left(X(z)+X(\rho z)+X(\rho^2 z)+\dots
+X(\rho^{N-1}z)\right)= \sum_n a_{Nk}z^{Nk} = \sum_n a_{Nk}(z^N)^k
= A(z^N)  .
$$
Here we used property $1+\rho+\rho^2+\dots \rho^{N-1}=0$,
which correctly for any root of degree $N$, $\rho \neq 1$, of unit.
Thus, in terms of power series wavelet expansion is
$$
A(z^N)=\frac{1}{N}\sum_{s=0}^{N-1}H_0(\rho^s z)X(\rho^s z), \quad
A_i(z^N)=\frac{1}{N}\sum_{s=0}^{N-1}H_i(\rho^s z)X(\rho^s z),
\quad i=1,2,\dots N-1.
$$
Reconstruction is made by other filters  $G_i(z)=\sum_n
\widetilde{g}_n^i z^n$, $i =0,1,2,\dots, N-1$, by formula
$$
a_n= \sum_k \widetilde{g}_{n-Nk}^0 a_{1,k}+\sum_{i=1}^{N-1}\sum_k
\widetilde{g}_{n-Nk}^i d_{1,k}^i .
$$
In terms of power series it means the following:
\begin{align*}
\sum_{i=0}^{N-1}G_i(z)A_i(z^N)=\frac{1}{N}\sum_{i=0}^{N-1}G_i(z)
\sum_{s=0}^{N-1}H_i(\rho^s z)X(\rho^s z)= \\*[-.5ex]
=\frac{1}{N}\sum_{s=0}^{N-1}\left(\sum_{i=0}^{N-1}G_i(z)H_i(\rho^s
z)\right)X(\rho^s z)=X(z).
\end{align*}
Therefore it is enough for perfect reconstruction, that it was
fulfilled equalities
$$
\sum_{i=0}^{N-1}G_i(z)H_i(z)=N \ \text{ at } s = 0, \quad
\sum_{i=0}^{N-1}G_i(z)H_i(\rho^s z)=N \ \text{ at } s = 1,2,\dots
N-1 .
$$
It means, that
$$
\small{\left( \begin{array}{ccc}
  G_0(z) &  \ldots & G_{N-1}(z) \\
  G_0(\rho z) &  \ldots & G_{N-1}(\rho z) \\
  \hdotsfor[2.5]{3} \\
  G_0(\rho^{N-1} z) &  \ldots & G_{N-1}(\rho^{N-1} z) \\
\end{array} \right)
\left(\begin{array}{cccc}
  H_0(z) & H_0(\rho z) & \ldots & H_0(\rho^{N-1} z) \\
  H_1(z) & H_1(\rho z) & \ldots & H_1(\rho^{N-1} z) \\
  \hdotsfor[2.5]{4} \\
  H_{N-1}(z) & H_{N-1}(\rho z) & \ldots & H_{N-1}(\rho^{N-1} z) \\
\end{array} \right) = N }.
$$

We have obtained the following fact.

{\bf Theorem 3.} {\it If a matrix of decomposition filters
$$
H(z)=\frac {1}{\sqrt{N}}\left(
\begin{array}{cccc}
  H_0(z) & H_0(\rho z) & \ldots & H_0(\rho^{N-1} z) \\
  H_1(z) & H_1(\rho z) & \ldots & H_1(\rho^{N-1} z) \\
  \hdotsfor[2.5]{4} \\
  H_{N-1}(z) & H_{N-1}(\rho z) & \ldots & H_{N-1}(\rho^{N-1} z) \\
\end{array} \right)  \eqno(16)
$$
is nonsingular, then exists the perfect reconstruction of a signal
by filters $G_i(z)$, $i = 0, 1, 2,\dots , N-1$, which matrix
$$
G(z)=\frac {1}{\sqrt{N}} \left(
\begin{array}{cccc}
  G_0(z) & G_0(\rho z) & \ldots & G_0(\rho^{N-1} z) \\
  G_1(z) & G_1(\rho z) & \ldots & G_1(\rho^{N-1} z) \\
  \hdotsfor[2.5]{4} \\
  G_{N-1}(z) & G_{N-1}(\rho z) & \ldots & G_{N-1}(\rho^{N-1} z) \\
\end{array} \right)  \eqno(17)
$$
is inverse transposed to matrix $H(z)$.}

{\bf Note.} In orthogonal case $H(z)$ is a unitary matrix, and
$G(z)$ is a complex conjugate matrix to $H(z)$.

{\bf 4. $B$-spline wavelets.} We shall assume, that $\varphi(x)$
is $B$-spline scaling function. As is known, it does not define
orthogonal multiresolution analysis. It is known also, that its
frequency function $H_0(z)$ is polynomial. We shall find wavelets
$\psi^1(x)\dots , \psi^{N-1}(x)$ and its polynomial frequency
functions, which lead to $N$-channel expansion of a input signal
and perfect reconstruction of the input signal by dual wavelets
$\widetilde{\psi}^1(x)\dots ,\widetilde{\psi}^{N-1}(x)$.

{\bf 4.1. General constructions.} Let $\varphi(x)$ is $B$-spline
scaling function and $H_0(z)$ is its frequency function. Let
$H_1(\omega),\dots , H_{N-1}(z)$ is frequency (polynomial)
functions of filters of expansion. According to the theorem 3, for
perfect reconstruction it is necessary, that matrix (16), where $z
= e^{-i\omega}$ è $\rho = e^{i2\pi/N}$ was nonsingular. We shall
include the factor $1/\sqrt{N}$ in frequency functions, as in the
formula (2). The matrix (16) has a special type. It is possible
become free of this special type of matrix $H(z)$ with used
Fourier transform on cyclical group $\mathbb{Z}/N\mathbb{Z} = \{1,
\rho , \rho^2,\dots , \rho^{N-1}\}$ \cite{BJ}. We shall define
$$
A_{i,j}(w)=\frac {1}{N}\sum_{z^N=w}z^{-j}H_i(z). \eqno(18)
$$
It is easy to check up, that the sum on the right depends from $w
= z^N$. The inverse transform is defined by the formula
$$
H_i(z)=\sum_{j=0}^{N-1}z^{j}A_{i,j}(z^N).
$$
Then the last relation can be represented as:
$$
H(z)=A(z^N)\left(
\begin{array}{cccc}
  1 & 1 & \ldots & 1 \\
  z & \rho z & \ldots & \rho^{N-1} z \\
  \hdotsfor[2.5]{4} \\
  z^{N-1} & \rho^{N-1} z^{N-1} & \ldots & \rho^{(N-1)^2} z^{N-1} \\
\end{array}
\right) = A(z^N)R(z). \eqno(19)
$$
Matrix $A(z^N)$ is already arbitrary nonsingular
matrix with polynomial elements. Now specific of the matrix $H(z)$ is
in the matrix $R(z)$. Setting matrix $A(z^N)$, we can construct special type
matrix $H(z)$ by (19) and also frequency functions $H_1(z),\dots ,
H_{N-1}(z)$ of wavelets, and hence wavelets
$\psi^1(x)\dots , \psi^{N-1}(x)$ themselves.
We shall assume, that polynomial frequency function $H_0(z)$ is defined, then the
first row of the matrix $A(z^N)$ is known,
$$
A_{0,j}(w)=\frac {1}{N}\sum_{z^N=w}z^{-j}H_0(z). \eqno(20)
$$
Now we construct remaining rows. As each element of first row of
the matrix $A(w)$ is a polynomial, the first row
$\alpha(w)=(A_{0,0}(w),\ \dots , A_{0,N-1}(w))$ can be represented
as
$$
\alpha(w)=\alpha_0 +\alpha_1 w +\dots
+\alpha_{g-1}w^{g-1},\eqno(21)
$$
where $\alpha_0, \alpha_1, \dots , \alpha_{g-1}$ there are $g$
vectors from $\mathbb{C}^N$.

{\bf Theorem 4.} {\it Let $\varphi(x)$ be a $B$-spline and
$H_0(z)$ -- its frequency function. If
$$
\alpha(w)=\alpha_0 +\alpha_1 w +\dots +\alpha_{g-1}w^{g-1}
$$
is expansion of the first row $\alpha(w)=(A_{0,0}(w),\dots ,
A_{0,N-1}(w))$, obtained by Fourier transform
(20) of frequency function $H_0(z)$, then
$$
\alpha_0 +\alpha_1 +\dots +\alpha_{g-1}= \frac 1N \cdot(1,1,\dots,
1).
$$
}

{\it Proof.} We shall use an induction with respect to degree $p$
of spline. If $p=0$, then $B$-spline is Haar function
$\varphi_0(x)$. Then
$$
H_0(z)= \frac {1}{N}\left(1+ z+ z^2+\dots + z^{N-1} \right).
$$
It easy to see that $A_{0,j}(w)=1/N$, therefore
$\alpha(w)=\alpha_0=N^{-1}\cdot(1,1,\dots, 1)$. Let's assume, that
the result has been proved for splines $\varphi_p(x)$ of a degree
$p$, and we shall prove it for $B$-splines of a degree $p+1$. As
is known, $\varphi_p(x)=\varphi_{p-1}(x)\ast\varphi_0(x)$. Hence
$\widehat{\varphi_p}(\omega)=(\widehat{\varphi_0}(\omega))^{p+1}$.
Therefore from refinement equation
$$
\widehat{\varphi_0}(\omega)=H_0\left( \frac{\omega}{N}\right)
\widehat{\varphi_0}\left(\frac{\omega}{N} \right)
$$
we get that the frequency function $H_{0,g}(z)$ of $B$-spline
$\varphi_g(x)$ of degree $g$ becomes:
$$
H_{0,p}(z)=(H_0(z))^{p+1}= \frac {1}{N^{p+1}}\left(1+ z+ z^2+\dots
+ z^{N-1} \right)^{p+1}.
$$
Let's present $H_{0,p}(z)$ as
$$
H_{0,p}(z)=\sum_{k=0}^{p}z^{kN}\sum_{i=0}^{N-1}a_{k,i}z^{i}=
\sum_{k=0}^{p}z^{kN}(a_{k,0}+a_{k,1}z+\dots + a_{k,N-1}z^{N-1}).
$$
Then $A_{0,j}(w)=\frac {1}{N}\sum_{z^N=w}z^{-j}H_{0,p}(z)$ is
sample of elements with degrees, multiply $z^{N}$ in polynomials
$H_{0,p}(z)$, $z^{-1}H_{0,p}(z)$, \dots , , $z^{-N+1}H_{0,p}(z)$.
Therefore
$$\begin{array}{lll}
  A_{0,0}= & \sum_{k=0}^{p}z^{kN}a_{k,0}= & a_{0,0}+a_{1,0}z^{N}+\dots+ a_{p,0}z^{pN}, \\
  A_{0,1}=  & \sum_{k=0}^{p}z^{kN}a_{k,1}= & a_{0,1}+a_{1,1}z^{N}+\dots+ a_{p,1}z^{pN}, \\
  \hdotsfor[2.5]{3} \\
  A_{0,N-1}=  & \sum_{k=0}^{p}z^{kN}a_{k,N-1}= & a_{0,N-1}+a_{1,N-1}z^{N}+\dots+ a_{p,N-1}z^{pN}. \\
\end{array}
$$
Hence,
$$\begin{array}{lll}
  \alpha_{0} & = & (a_{0,0}, a_{0,1},\dots, a_{0,N-1}), \\
  \alpha_{1} & = & (a_{1,0}, a_{1,1},\dots, a_{1,N-1}), \\
  \hdotsfor[2.5]{3} \\
  \alpha_{p}  & = & (a_{p,0}, a_{p,1},\dots, a_{p,N-1}). \\
\end{array}\eqno(22)
$$
From the induction hypothesis, the sum of these vectors is a
vector $(1,1,\dots, 1)/N$.

Now let's consider $B$-spline $\varphi_{p+1}(x)$ of degree $p+1$.
Its frequency function $H_{0,p+1}(z)$ building from $H_{0,p}(z)$
by multiplication to a polynomial $H_{0}(z)$,
$$
H_{0,p+1}(z)=\frac 1N H_{0,p}(z)(1+ z+ z^2+\dots +
z^{N-1}).\eqno(23)
$$
Factor $1/N$ in this expression we shall not take into account
yet. As well as earlier the Fourier transform $A_{0,j}(w)=\frac
{1}{N}\sum_{z^N=w}z^{-j}H_{0,p+1}(z)$ is sample of elements with
degrees, multiply $z^{N}$ in polynomials $H_{0,p+1}(z)$,
$z^{-1}H_{0,p+1}(z)$, \dots , $z^{-N+1}H_{0,p+1}(z)$. We shall
uncover brackets in (23) (without factor $1/N $) and make such
sample separately for each part.

For the factor $z^0$ we use Fourier transform to $H_{0,p}(z)$.
Then $A_{0,j}(w)=\frac {1}{N}\sum_{z^N=w}z^{-j}H_{0,p}(z)$ is
sample of elements with degrees, multiply $z^{N}$ in polynomials
$H_{0,p}(z)$, $z^{-1}H_{0,p}(z)$, \dots , $z^{-N+1}H_{0,p}(z)$.
Corresponding vector coefficients coincide with found earlier in
(22), $\alpha_{0,0}=\alpha_{0}$, $\alpha_{0,1}=\alpha_{1}$, \dots,
$\alpha_{0,p}=\alpha_{p}$. From the induction hypothesis, the sum
of these vectors is equal $(1,1,\dots, 1)/N$.

For the factor $z^1$ we use Fourier transform to $H_{0,p}(z)z$.
Then $A_{0,j}(w)=\frac {1}{N}\sum_{z^N=w}z^{-j}H_{0,p}(z)z=\frac
{1}{N}\sum_{z^N=w}z^{-j+1}H_{0,p}(z)$ is sample of elements with
degrees, multiply $z^{N}$ in polynomials $H_{0,p}(z)z$,
$H_{0,p}(z)$, $z^{-1}H_{0,p}(z)$ \dots , $z^{-N+2}H_{0,p}(z)$.
Let's consider more in detail sample of elements with degrees,
multiple $z^{N}$ in the element $H_{0,p}(z)z$. We have,
$$
H_{0,p}(z)z=\sum_{k=0}^{p}z^{kN}\sum_{i=0}^{N-1}a_{k,i}z^{i}\cdot
z =\sum_{k=0}^{p}z^{kN}(a_{k,0}z+a_{k,1}z^{2}+\dots+
a_{k,N-2}z^{N-1}+a_{k,N-1}z^{N}).
$$
Then
$$
\begin{array}{llllllll}
  A_{0,0} & = & 0 & + a_{0,N-1}z^{N} & + a_{1,N-1}z^{2N} & +\dots  & + a_{p-1,N-1}z^{pN}+ a_{p,N-1}z^{(p+1)N}, \\
  A_{0,1} & =  & a_{0,0} & +a_{1,0}z^{N} & +a_{2,0}z^{2N} & +\dots  & + a_{p,0}z^{pN}, \\
  \hdotsfor[2.5]{6} \\
  A_{0,N-1} & = & a_{0,N-2} & +a_{1,N-2}z^{N} & +a_{2,N-2}z^{2N} & +\dots  & + a_{p,N-2}z^{pN}. \\
\end{array}
$$
Hence,
$$\begin{array}{lllcccccr}
  \alpha_{1,0} &  = & ( & 0, & a_{0,0}, & a_{0,1}, & \dots , & a_{0,N-2} &), \\
  \alpha_{1,1} &  = & ( & a_{0,N-1}, & a_{1,0}, & a_{1,1}, & \dots , & a_{1,N-2}&), \\
  \hdotsfor[2.5]{9} \\
  \alpha_{1,p} &  = & ( & a_{p-1,N-1}, & a_{p,0}, & a_{p,1} & \dots, & a_{p,N-2}&), \\
  \alpha_{1,p+1} &  = & ( & a_{p,N-1}, & 0, & 0, & \dots , & 0&). \\
\end{array}
$$
The sum of these vectors is the vector $(1,1,\dots, 1)/N$ again.

Let's continue this procedure. For the factor $z^{N-1}$ is used
Fourier transform to $H_{N-1,0,p+1}(z)=H_{0,p}(z)z^{N-1}$. Then
$A_{1,0,j}(w)=\frac {1}{N}\sum_{z^N=w}z^{-j+N-1}H_{0,p}(z)$ is
sample of elements with degrees, multiple $z^{N}$ in polynomials
$H_{0,p}(z)z^{N-1}$, $H_{0,p}(z)z^{N-2}$, \dots , $H_{0,p}(z)z$
$H_{0,p}(z)$. In this case vector coefficients are
$$\begin{array}{lllcccccr}
  \alpha_{N-1,0} &  = & ( & 0, & 0, & \dots , & 0,  & a_{0,0} &), \\
  \alpha_{N-1,1} &  = & ( & a_{0,1}, & a_{0,2}, & \dots , & a_{0,N-1}, & a_{1,0}&), \\
  \hdotsfor[2.5]{9} \\
  \alpha_{N-1,p} &  = & ( & a_{p-1,1}, & a_{p-1,2}, & \dots, & a_{p-1,N-1}, & a_{p,0}&), \\
  \alpha_{N-1,p+1} &  = & ( & a_{p,1}, & a_{p,2}, & \dots , & a_{p,N-1},& 0, &). \\
\end{array}
$$
The sum of these vectors is the vector $(1,1,\dots, 1)/N$ again.

Summarizing all obtained vector coefficients of their $N$ groups,
we get, that their sum is a vector $(1,1,\dots, 1)$. Taking into
account the factor $1/N$ omitted earlier, we obtain, that the sum
of all vector coefficients is equal $(1,1,\dots, 1)/N$.

The theorem is proved for $B$-spline with support on interval
$[0,p+1]$. In a case when the support of $B$-spline is centered,
as in section 1.3, the theorem is proved similarly.

Let's return to expansion (21) of first row
$\alpha(w)=(A_{0,0}(w),\dots , A_{0,N-1}(w))$. It can be
complemented up to a nonsingular matrix with the help simple
procedure which is generalization of a similar construction in an
orthogonal case \cite{BJ}. As vectors $\alpha_i$ in (21) can be
linearly dependent we shall consider more general case,
$$
\alpha(w)=\alpha_0 p_0(w) +\alpha_1 p_1(w) +\dots +\alpha_{g-1}
p_{g-1}(w),
$$
where $\alpha_0, \alpha_1, \dots , \alpha_{g-1}$ are $g$ linearly
independent vectors from $\mathbb{C}^N$ è $p_0(w), p_1(w), \dots ,
p_{g-1}(w)$ are some polynomials.

{\bf Theorem 5.} {\it Let $\alpha= (\alpha_0, \alpha_1, \dots ,
\alpha_{g-1})$ are $g$ linearly independent vectors from
$\mathbb{C}^N$ and $p_0(w), p_1(w), \dots , p_{g-1}(w)$ are some
polynomials, any of which not vanish on the unit circle $w =
e^{it}$. Then there is polynomial loop $A(w)$ on group
$GL(N,\mathbb{C})$ such, that the first row of $A(w)$ is
$\sum_{i=0}^{g-1}\alpha_i p_{i}(w)$. Thus the degree of loop
$A(w)$ is at most the maximum degree of polynomials $p_i(w)$.}

{\it Proof.} We shall use induction with respect to maximum degree
$k$ of polynomials $p_i(w)$. Without loss of generality it is
possible to assume, that the sum $\sum_{i=0}^{g-1}\alpha_i
p_{i}(w)$ is ordered on increase of degrees of polynomials
$p_i(w)$. If $k = 0$ then $\alpha_0$ is some nonzero row vector,
and we can put $A(w) = A$, where $A$ is a nonsingular matrix with
first row $\alpha_0$.

Let's assume, that $k >0$ and that the result has been proved for
all smaller $k$. We can suppose, that $\alpha_{g-1}\neq 0$. Let
$P$ is an one-dimensional projection to vector $\alpha_{g-1}$
along a subspace formed by the first vectors $\alpha_0, \alpha_1,
\dots , \alpha_{g-2}$. Then
$$
\alpha_{g-1} P = \alpha_{g-1},\quad \alpha_{i}P = 0,\quad
i=0,1,\dots g-2.
$$
Let's remind, that the projection operator $P$ has the following
properties: $ P^2 = P$, \ $P(1-P) = 0$,\ $(1-P)^2 = 1-P$. Now we
shall define
\begin{equation*}
\begin{split}
  \beta(w) & =\alpha(w)(1-P+p_{g-1}(w)^{-1}P)= \\
    & = \left(\alpha_0 p_0(w) +\alpha_1 p_1(w) +\dots +\alpha_{g-1}
p_{g-1}(w) \right)(1-P+p_{g-1}(w)^{-1}P)= \\
    & = p_{g-1}(w)^{-1}\alpha_0 P + p_0(w)\alpha_0(1-P)+\\
    & +p_1(w)\alpha_1(1-P)+p_1(w)p_{g-1}(w)^{-1}\alpha_1
P+\dots +p_{g-1}(w)p_{g-1}(w)^{-1}\alpha_{g-1} P +  \\
    & +p_{g-1}(w)\alpha_{g-1}(1- P)= \\
    & =\alpha_{g-1}+\alpha_0 p_0(w) +\alpha_1 p_1(w) +\dots
+\alpha_{g-2} p_{g-2}(w) .
\end{split}
\end{equation*}
Therefore $\beta(w)$ is a polynomial on $w$ of degrees, smaller
$k$, as $\deg(p_{g-2}(w))<\deg(p_{g-1}(w))=k$. By induction
hypothesis, exists polynomial loop $B(w)$ of degrees
$\deg(p_{g-2}(w))$ such, that first row $B(w)$ is $\beta(w)$.
Then, supposing $A(w) = B(w)(1-P + p_{g-1}(w)P)$, from the
obtained expression for $\beta(w)$ follows, that first row $A(w)$
is $\alpha(w)$. It complete the proof of the Theorem 5.

{\bf 4.2. $B$-spline filters of expansion and reconstruction.} As
is known, frequency function for $B$-spline scaling function is
polynomial and does not define orthogonal multi\-resolu\-tion
analysis. Constructions mentioned above allow to build
nonorthogonal wavelets with polynomial frequency functions and to
find dual functions which ensure perfect reconstruction of a
signal. If polynomial frequency function $H_0 (z) $ is set, it is
possible to be considered, that the first row of a matrix $A (z^N)
$ is known
$$
A_{0,j}(w)=\frac {1}{N}\sum_{z^N=w}z^{-j}H_0(z).
$$
Remaining rows can be constructed as in the proof of the theorem
5. As each element of the first row of $A(w)$ is a polynomial, the
first row $\alpha(w)=(A_{0,0}(w),\ \dots , A_{0,N-1}(w))$ can be
represented as
$$
\alpha(w)=\alpha_0 +\alpha_1 w +\dots +\alpha_{g-1}w^{g-1},
$$
where $\alpha_0, \alpha_1, \dots , \alpha_{g-1}$ are $g$ vectors
from $\mathbb{C}^N$. Let $P_{1}$ is an one-dimensional projection
to a vector $\alpha_{g-1}$ along subspace formed by the first
vectors $\alpha_0, \alpha_1, \dots , \alpha_{g-2}$ Then
\begin{equation*}
\begin{split}
  \beta(w) & =\alpha(w)(1-P_{1}+w^{-g+1}P_{1})= \\
    & = \left(\alpha_0 +\alpha_1 w +\dots +\alpha_{g-1}w^{g-1} \right)(1-P_{1}+w^{-g+1}P_{1})= \\
    & =w^{-g+1}\alpha_0 P_{1} + \alpha_0(1-P_{1})+\\
    & +w^{1}\alpha_1(1-P_{1})+w^{1}w^{-g+1}\alpha_1
P_{1} +\dots +w^{g-1} w^{-g+1}\alpha_{g-1} P_{1} +  \\
    & +w^{g-1}\alpha_{g-1}(1- P_{1})= \\
    & =\alpha_0+\alpha_{g-1} +\alpha_1 w +\dots
+\alpha_{g-2}w^{g-2} .
\end{split}
\end{equation*}
Then
$$
\alpha(w)=(\alpha_0 +\alpha_{g-1} +\alpha_1 w +\dots +\alpha_{g-2}
w^{g-2})(1-P_{1}+w^{g-1}P_{1}).
$$
Continuing this procedure, we obtain the following expansion:
$$
\alpha(w)= (\alpha_0 +\alpha_1 +\dots
+\alpha_{g-1})\prod_{k=1}^{g-1}(1-P_{k}+w^{g-k}P_{k}),\eqno(24)
$$
where $P_{k}$ is an one-dimensional projection to a vector $\alpha_{g-k}$
along a subspace formed by remaining vectors $\alpha_i$.
In the theorem 4 it is shown, that the vector $\alpha_0 +\alpha_1 +\dots
+\alpha_{g-1}$ is $(1,1,\dots, 1)/N$. This vector needs to be complemented up
to a nonsingular matrix $A_0$ of the order $N$.
We shall specify two most simple methods:
$$
A_0=\left(%
\begin{array}{cccccc}
  1 & 1 & 1 & \ldots & 1 & 1 \\
  0 & 1 & 0 & \ldots & 0 & 0 \\
  0 & 0 & 1 & \ldots & 0 & 0 \\
  \hdotsfor[2.5]{6} \\
  0 & 0 & 0 & \ldots & 0 & 1 \\
\end{array}%
\right)\quad  \text { è } \quad
A_0=\left(%
\begin{array}{cccccc}
  1 & 1 & 1 & \ldots & 1 & 1 \\
  1 & -1 & 0 & \ldots & 0 & 0 \\
  1 & 1 & -2 & \ldots & 0 & 0 \\
  \hdotsfor[2.5]{6} \\
  1 & 1 & 1 & \ldots & 1 & -N+1 \\
\end{array}%
\right).
$$
In the second case the matrix $A_0$ has orthogonal rows and its
possible to make an orthogonal by normalization of rows. Thus
there are many such matrixes $A_0 $. Each such matrix can be
obtained from some one by multiplication at the left on an
arbitrary nonsingular matrix $g$ of the form
$$
g=\left(%
\begin{array}{cccc}
  1 & 0 &  \ldots  & 0 \\
  a_{1,0} & a_{1,1} & \ldots  & a_{1,N-1} \\
  \hdotsfor[2.5]{4} \\
  a_{N-1,0} & a_{N-1,1} & \ldots  & a_{N-1,N-1} \\
\end{array}%
\right).
$$
After a choice of constant matrix $A_0$, the matrix $A(z^N)$ is
given by formula
$$
A(w)= \frac 1N
A_0\prod_{k=1}^{g-1}(1-P_{k}+w^{g-k}P_{k}),\eqno(25)
$$
where $w=z^N$. Now the matrix of initial filters $H(z)$ become $H(z)=A(z^N)R(z)$,
where $R(z)$ is the complementary matrix indicated in the formula (19).
The matrix $G(z)$ of reconstruction filters is inverse transposed
to a matrix initial filters. As $H(z)=A(z^N)R(z)$, then
$$
G^t(z)=H^{-1}(z)=R^{-1}A^{-1}=\frac {1}{N} \left(
\begin{array}{rrrr}
  1 & z^{-1} & \ldots & z^{-(N-1)} \\
  1 & \rho^{-1}z^{-1} & \ldots & \rho^{-(N-1)}z^{-(N-1)}\\
  \hdotsfor[2.5]{4} \\
  1 & \rho^{-(N-1)} z^{-1} & \ldots & \rho^{-(N-1)^2} z^{-(N-1)} \\
\end{array}
\right)A(z^N)^{-1} .
$$
Last formula enables to find filters of reconstruction (dual wavelets).
We shall show it on examples of $B$-splines of order 1 and 2 with
parameter of scaling $N=3$.

{\bf Example 4.2.1.}
Let $N=3$. We shall consider $B$-spline of order 1. As we notes earlier,
frequency function looks like
$$
H_0(z)=\frac 19z^{-2}(1+2z+3z^2 +2z^3 +z^4 ).
$$
Let's calculate
$$
A_{0,j}(w)=\frac {1}{N}\sum_{z^N=w}z^{-j}H_0(z),\qquad w\in
\mathbb{T}, \quad j=0,1,2\ .
$$
Let's put $w=e^{-i\omega}$ and $\rho=e^{i2\pi/3}$ then if $z_0=e^{-i\omega/3}$
summation is made on the following values $z\in\{ z_0, \rho z_0, \rho^2 z_0 \}$.
In the further evaluations we shall omit a factor $1/9$. Then we shall put,
$$
m_0(z)=z^{-2}(1+2z+3z^2 +2z^3 +z^4 )= z^{-2}+2z^{-1}+3 +2z +z^2 .
$$
also we shall calculate a Fourier transform on cyclical group of
the third order  $\{ 1, \rho , \rho^2 \}$,
$$
A_{0,0}(w)=\frac {1}{3}\sum_{z^3=w}m_0(z) =\frac 13 9 =3,\qquad
A_{0,1}(w)=\frac {1}{3}\sum_{z^3=w}z^{-1}m_0(z) =\frac 13
(3z^{-3}+6) = z^{-3}+2,
$$
$$
A_{0,2}(w)=\frac {1}{3}\sum_{z^3=w}z^{-2}m_0(z)= \frac
{1}{3}(6z^{-3}+3)= 2z^{-3}+1\ .
$$
We decompose the first row on degrees with vector coefficient,
$$
\left(A_{0,0}(w),A_{0,1}(w),A_{0,2}(w)\right)=
\left(3,2+z^{-3},1+2z^{-3}\right)= (3,2,1)+(0,1,2)w^{-1},
$$
$$
\alpha_0=(3,2,1), \quad  \alpha_1=(0,1,2), \quad z^{-3}=w^{-1}\ .
$$
Let $P$ there is an one-dimensional projection to a vector $\alpha_1$
along $\alpha_0$ i.e., $\alpha_1 P = \alpha_1$ è $\alpha_0 P = 0$.
Now we shall define
$$
\beta(w) = \alpha(w)(1-P + wP)= (\alpha_0 +\alpha_1 w^{-1})(1-P +
wP)=\alpha_0 +\alpha_1.
$$
It is obvious, that $(1-P + w^{-1}P)$ is a inverse matrix for $(1-P + wP)$, therefore
$$
\alpha(w)= \beta(w)(1-P + w^{-1}P)= (\alpha_0 +\alpha_1)(1-P +
w^{-1}P), \text{ where } \alpha_0 +\alpha_1= (3,3,3) .
$$
We complement this first row $\alpha_0 +\alpha_1$ up to a
nonsingular matrix $A_0$ for example, unit rows $(0, 1, 0)$ and
$(0, 0, 1)$. Then $A(w)=A_0(1-P + w^{-1}P)$. We shall find
expression $(1-P + w^{-1}P)$ of an operator. For this purpose we
shall add vectors $\alpha_0=(3,2,1)$ and $\alpha_1=(0, 1, 2)$ up
to basis $\mathbb{C}^3$. Let $\alpha_2=(0,0,1)$. In this basis the
operator $(1-P + w^{-1}P)$ becomes
$$
1-P + w^{-1}P= \left(%
\begin{array}{ccc}
  1 & 0 & 0 \\
  0 & w^{-1} & 0 \\
  0 & 0 & 1 \\
\end{array}%
\right).
$$
Let's find expression in standard basis $e_0,e_1,e_2$ of spaces
$\mathbb{C}^3$. For this purpose we shall find a transition matrix
$B$ from one basis, to another $\alpha = Be$, $e=B^{-1}\alpha$,
$$
B=\left(%
\begin{array}{ccc}
  3 & 2 & 1 \\
  0 & 1 & 2 \\
  0 & 0 & 1 \\
\end{array}%
\right),
\quad B^{-1}=\left(%
\begin{array}{rrr}
  1/3 & -2/3 & 1 \\
  0 & 1 & -2 \\
  0 & 0 & 1 \\
\end{array}%
\right) .
$$
Let $P_{\alpha k}^j$ is the matrix of operator in basis
$(\alpha_0,\alpha_1, \alpha_2)$. The operator matrix in standard
basis is evaluated with help of the formula $P_e=B^{-1}P_\alpha
B$, $P_{e s}^q= (B^{-1})_{s}^k P_{\alpha k}^j B_j^q$. Then the
operator $(1-P + w^{-1}P)$ in standard basis has the matrix
$$
P_e(w)= \left(\begin{array}{ccc}
  1 & 2/3 (1-w^{-1}) & 4/3 (1-w^{-1}) \\
  0 & w^{-1} & -2(1-w^{-1}) \\
  0 & 0 & 1 \\
\end{array}\right) .
$$
It allows to calculate polynomial matrix $A(w)$:
$$
A(w)=A_0((1-P + w^{-1}P))=A_0P_e(w)= \left(\begin{array}{ccc}
  3 & 2+w^{-1} & 1+2w^{-1} \\
  0 & w^{-1} & -2(1-w^{-1}) \\
  0 & 0 & 1 \\
\end{array}\right) .
$$
Let's calculate the matrix of frequency functions by the formula $H(z)=A(z^N)R(z)$,
$$
\left( \begin{array}{ccc}
  m_0(z) & m_0(\rho z) & m_0(\rho^2 z) \\
  m_1(z) & m_1(\rho z) & m_1(\rho^2 z)\\
  m_2(z) & m_2(\rho z) & m_2(\rho^2 z) \\
\end{array} \right)=
\left(\begin{array}{ccc}
  3 & 2+z^{-3} & 1+2z^{-3} \\
  0 & z^{-3} & -2(1-z^{-3}) \\
  0 & 0 & 1 \\
\end{array}\right)
\left(\begin{array}{ccc}
  1 & 1 & 1 \\
  z & \rho z & \rho^2 z \\
  z^2 & \rho^2 z^2 & \rho^4 z^2 \\
\end{array}\right) .
$$
We get, $m_0(z)=z^{-2}+2z^{-1}+3+2z+z^{2}$,
$m_1(z)=z^{-2}+2z^{-1}-2z^{2}$, $m_2(z)=z^{2}$. Therefor
$$
H_0(z)=\frac 19\left(z^{-2}+2z^{-1}+3+2z+z^{2}\right), \quad
H_1(z)=\frac 19\left(z^{-2}+2z^{-1}-2z^{2}\right), \quad
H_2(z)=\frac 19z^{2} .
$$
The matrix $G(z)$ of reconstruction filters is an inversion of
matrix $H(z)$, $G(z)=R(z)^{-1}A(z^n)^{-1}$,
$$
\small{\left( \begin{array}{rrr}
  G_0(z) & G_1(z) & G_2(z) \\
  G_0(\rho z) & G_1(\rho z) & G_2(\rho z) \\
  G_0(\rho^2 z) & G_1(\rho^2 z) & G_2(\rho^2 z) \\
\end{array} \right){=}
3\!\left(\begin{array}{rrr}
  1 & z^{-1} & z^{-2} \\
  1 & \rho^{-1} z^{-1} & \rho^{-2} z^{-2} \\
  1 & \rho^{-2} z^{-1} & \rho^{-4} z^{-2} \\
\end{array}\right)\!
 \left(\begin{array}{ccc}
  1/3 & -(1{+}2z^3)/3 & (1{-}4z^3)/3 \\
  0 & z^3 & -2(1{-}z^3) \\
  0 & 0 & 1 \\
\end{array}\right) .}
$$
Thus it is necessary to take into account, that in wavelet expansion take part
the conjugate filters of wavelets. It means, that we must to convert a matrix $H(z)$,
with a complex conjugate elements. As filter coefficient is real,
it is enough to make substitution $z=e^{-i\omega}$ on $z^{-1}$. We choose elements
of the first row and obtain the following filters of reconstruction,
$$
G_0(z) = 1,\qquad  G_1(z) = - 2z^{-3}+ 3z^{-2}-1 ,\qquad G_2(z) =
-4z^{-3} + 6z^{-2} + 1 - 6z+3z^2\ .
$$
Thus, reconstruction filters have the following nonzero elements:
$$
\begin{array}{l}
  \widetilde{h}_0=\sqrt{3}, \\
  \widetilde{g}_{-3}^1=-2\sqrt{3},\quad \widetilde{g}_{-2}^1=3\sqrt{3},\quad \widetilde{g}_{0}^1
  =-\sqrt{3},  \\
  \widetilde{g}_{-3}^2= -4\sqrt{3},\quad \widetilde{g}_{-2}^2=6\sqrt{3},\quad \widetilde{g}_{0}^2
  =\sqrt{3}, \quad \widetilde{g}_{1}^2= -6\sqrt{3},\quad \widetilde{g}_{2}^2= 3\sqrt{3}\ . \\
\end{array}
$$

{\bf Example 4.2.2.}
Let $N=3$. We shall consider $B$-spline of order 2. Frequency function of
$\varphi_2(x)$ is
$$
H_0(z)=\frac {1}{3^3}z^{-2}\left(1+z+z^2\right)^3= \frac
{1}{27}\left(z^{-2}+ 3z^{-1}+6+7z+6z^2+3z^3+z^4\right).
$$

Let's calculate
$$
A_{0,j}(w)=\frac {1}{N}\sum_{z^N=w}z^{-j}H_0(z),\qquad w\in
\mathbb{T}, \quad j=0,1,2\ .
$$
Let's put $w=e^{-i\omega}$ and $\rho=e^{i2\pi/3}$, then if
$z_0=e^{-i\omega/3}$ summation is made on the following values
$z\in\{ z_0, \rho z_0, \rho^2 z_0 \}$. In the further evaluations
we shall omit a factor $1/27$. Then we shall put,
$$
m_0(z)=z^{-2}+ 3z^{-1}+6+7z+6z^2+3z^3+z^4
$$
also we shall calculate a Fourier transform on cyclical group of
the third order $\{ 1, \rho , \rho^2 \}$,
$$
A_{0,0}(w)=\frac {1}{3}\sum_{z^3=w}m_0(z) =6+3z^3,\qquad
A_{0,1}(w)=\frac {1}{3}\sum_{z^3=w}z^{-1}m_0(z) =z^{-3}+7+z^3,
$$
$$
A_{0,2}(w)=\frac {1}{3}\sum_{z^N=w}z^{-2}m_0(z)= 3z^{-3}+6\ .
$$
We decompose the first row,
$$
\left(A_{0,0}(w),A_{0,1}(w),A_{0,2}(w)\right)=
\left(6+3w,w^{-1}+7+w,3w^{-1}+6\right),
$$
$$
\alpha_{-1}=(0,1,3), \quad  \alpha_0=(6,7,6), \quad
\alpha_1=(3,1,0), \quad z^3=w\ .
$$
Then
$$
\alpha (w)= \alpha_{-1}w^{-1}+ \alpha_0+ \alpha_1 w.
$$
Let $P_1$ is an one-dimensional projection to the vector $\alpha_1$ along
$\alpha_{-1}$ and $\alpha_0$ i.e. $\alpha_1 P_1 = \alpha_1$  and
$\alpha_{-1} P_1 = 0$, $\alpha_0 P_1 = 0$. Now we shall define
$$
\beta(z) = \alpha(z)(1-P_1 + w^{-1}P_1)= (\alpha_{-1}w^{-1}+
\alpha_0+ \alpha_1 w)(1-P_1 +
w^{-1}P_1)=\alpha_{-1}w^{-1}+\alpha_0 +\alpha_1.
$$
Therefore
$$
\alpha(z)= (\alpha_{-1}w^{-1}+\alpha_0 +\alpha_1)(1-P_1 + w P_1).
$$
Now let $P_2$ is an one-dimensional projection to a vector $\alpha_{-1}$
along $\alpha_0+\alpha_{1}$. Then
$$
\alpha(z)= (\alpha_{-1}+\alpha_0 +\alpha_1)(1-P_2 + w^{-1}
P_2)(1-P_1 + w P_1).
$$
We complement the first row $\alpha_{-1}+\alpha_0 +\alpha_1
=(9,9,9)$ up to a nonsingular matrix $A_0$, for example, by unit
rows $(0, 1, 0)$ and $(0, 0, 1)$. Then $A(w)=A_0(1-P_2 + w^{-1}
P_2)(1-P_1 + w P_1)$. Let's find expressions of projective
operators in standard basis $\mathbb{C}^3$. First we shall
consider basis $\mathbb{C}^3$, consisting of vectors
$\alpha_{-1}$, $\alpha_0$, $\alpha_1$. In this base operators of
projection $(1-P_1 + wP_1)$ become
$$
1-P_1 + wP_1= \left(%
\begin{array}{ccc}
  1 & 0 & 0 \\
  0 & 1 & 0 \\
  0 & 0 & w \\
\end{array}%
\right), \quad
1-P_2 + w^{-1}P_2= \left(%
\begin{array}{ccc}
  w^{-1} & 0 & 0 \\
  0 & 1 & 0 \\
  0 & 0 & 1 \\
\end{array}%
\right).
$$
Let's find expression of projective operators in standard basis $e_0,e_1,e_2$ of spaces $\mathbb{C}^3$.
For this purpose we shall find a transition matrix $B$ from one basis, to another
$\alpha = Be$, $e=B^{-1}\alpha$.
$$
B=\left(%
\begin{array}{ccc}
  0 & 1 & 3 \\
  6 & 7 & 6 \\
  3 & 1 & 0 \\
\end{array}%
\right), \quad B^{-1}= \left( \begin {array}{rrr}
2/9&-1/9&5/9\\
-2/3&1/3&-2/3\\
5/9&-1/9&2/9
\end{array} \right).
$$
Let $Q_{\alpha k}^j$ is the matrix of operator in basis
$(\alpha_0, \alpha_1, \alpha_2)$. The matrix of an operator in
standard basis is evaluated with help of the formula
$Q_e=B^{-1}Q_\alpha B$, $Q_{e s}^q= (B^{-1})_{s}^k Q_{\alpha k}^j
B_j^q$. Then operators $1-P_1 + wP_1$ $1-P_2 + w^{-1}P_2$ in
standard base have matrixes
$$
\left( \begin {array}{ccc} -2/3+5/3\,w&-5/9+5/9\,w&0\\
2-2\,w&5/3-2/3\,w&0\\
-2/3+2/3\,w &-2/9+2/9\,w&1\end {array} \right),\quad
\left( \begin {array}{ccc} 1&2/9\,{w}^{-1}-2/9&2/3\,{w}^{-1}-2/3\\
0&-2/3\,{w}^{-1}+5/3&-2\,{w}^{-1}+2\\
0&5/9\,{w}^{-1}-5/9&5/3\,{w}^{-1}-2/3\end{array}
 \right)
.
$$
It allows to calculate the polynomial matrix $A(w)=A_0(1-P_1 +
wP_1)(1-P_2 + w^{-1}P_2)$:
$$
A(w) = \left(\begin {array}{ccc} 6+3w & (1+{w}^{2}+7w)/w & 3(1+2w)/w\\
2-2\,w & -(2-7\,w+2\,w^2)/(3w) & 2\,(-1+w)/w\\
(-2+2\,w)/3 & (-1+w)(-5+2\,w)/(9w) &(5-2\,w)/(3w)\end
{array}\right) .
$$
We compute the matrix of frequency functions by the formula
$H(z)=A(z^N)R(z)$. The first column gives us frequency functions
of wavelets  $\psi^1(x)$ and $\psi^2(x)$:
$$
H_0(z)=\frac {1}{27z^2}\left(1 +3\,z +6\,{z}^{2}
+7\,{z}^{3}+6\,{z}^{4} +3\,{z}^{5} +{z}^{6}\right),
$$
$$
H_1(z)=\frac {1}{81z^2}\left(-2 -6\,z +6\,{z}^{2} +7\,{z}^{3}
+6\,{z}^{4} -6\,{z}^{5} -2\,{z}^{6}\right),
$$
$$
H_2(z)=\frac {1}{243 z^2}\left(5 +15\,z -6\,{z}^{2} -7\,{z}^{3}
-6\,{z}^{4} +6\,{z}^{5} +2\,{z}^{6}\right).
$$
The matrix $G(z)$ of filters of reconstruction is an inversion
of matrix $H(z)$, $G(z)=R(z)^{-1}A(z^n)^{-1}$. Thus it is necessary
to take into account, that in wavelet expansion take part the conjugate
filters of wavelets. It means, that we must to convert the matrix $H(z)$,
with a complex conjugate elements. As the filter coefficient is real,
it is enough to make substitution $z=e^{-i\omega}$ on $z^{-1}$. We choose elements
of the first row and obtain the following filters of reconstruction:
$$
G_0(z)= \frac
{z^5}{3}\left(-2\,z^{-5}+6\,z^{-4}-2\,z^{-3}+5\,z^{-2}-6\,z^{-1}+2
\right),
$$
$$
G_1(z)= z^5
\left(2\,z^{-8}-6\,z^{-7}+5\,z^{-6}-z^{-5}+3\,z^{-4}-z^{-3}-10\,z^{-2}+
12\,z^{-1} -4 \right),
$$
$$
G_2(z)= z^5
\left(6\,z^{-8}-18\,z^{-7}+15\,z^{-6}-15\,z^{-2}+18\,z^{-1}-6\right).
$$

{\bf Note.} The obtained expressions depend on a choice of a
matrix $A_0$. For example, if to take a matrix $A_0$ with
orthogonal rows,
$$
A_0=\left(\begin{array}{ccc}
  9 & 9 & 9 \\
  1 & -1 & 0 \\
  1 & 1 & -2 \\
\end{array} \right),
$$
then
$$
H_1(z)=\frac {1}{243z^2}\left(8 +24\,z -24\,{z}^{2} -28\,{z}^{3}
-24\,{z}^{4} +33\,{z}^{5} +11\,{z}^{6}\right),
$$
$$
H_2(z)=\frac {1}{243 z^2}\left(-14 -42\,z +24\,{z}^{2}
+28\,{z}^{3} +24\,{z}^{4} -15\,{z}^{5} -5\,{z}^{6}\right).
$$
Reconstruction filters is
$$
G_0(z)= \frac
{z^5}{27}\left(8\,{z}^{-8}-24\,{z}^{-7}+20\,{z}^{-6}-19\,{z}^{-5}+57\,{z}^{-4}-
19\,{z}^{-3} +20\,{z}^{-2}-24\,z^{-1}+8 \right),
$$
$$
G_1(z)= \frac {z^5}{2}
\left(-2\,{z}^{-8}+6\,{z}^{-7}-5\,{z}^{-6}+{z}^{-5}-3\,{z}^{-4}+{z}^{-3}+
10\,{z}^{-2}-12\,z^{-1}+4\right),
$$
$$
G_2(z)= \frac {z^5}{6}
\left(-10\,{z}^{-8}+30\,{z}^{-7}-25\,{z}^{-6}-{z}^{-5}+3\,{z}^{-4}-{z}^{-3}+20\,{z
}^{-2}-24\,z^{-1}+8\right).
$$

\end{document}